\documentclass[11pt,thmsa,emstex]{article}
\usepackage[left=3cm,right=3cm,top=3cm,bottom=3cm]{geometry}
\usepackage{amsfonts}
\usepackage{t1enc}
\usepackage{amssymb}
\usepackage{epsfig}
\usepackage[french,english]{babel}
\usepackage{amsmath}
\usepackage{graphics}
\usepackage{layout}
\usepackage{latexsym}
\usepackage{color}

\newtheorem{theorem}{Theorem}[section]
\newtheorem{lemma}[theorem]{Lemma}
\newtheorem{corollary}[theorem]{Corollary}
\newtheorem{proposition}[theorem]{Proposition}
\newtheorem{example}[theorem]{Example}
\newtheorem{remark}[theorem]{Remark}

\newtheorem{hypothesis}[theorem]{Hypothesis}

\newtheorem{exercice}[theorem]{Exercice}
\newcommand*{\QEDB}{\hfill\ensuremath{\square}}

\def\bit{\begin{itemize}}

\def\eit{\end{itemize}}

\reversemarginpar   

\def\bc{\begin{center}}

\def\ec{\end{center}}

\def\bthm{\begin{theorem}}

\def\ethm{\end{theorem}}

\def\bcor{\begin{corollary}}

\def\ecor{\end{corollary}}

\def\bprop{\begin{proposition}}

\def\eprop{\end{proposition}}

\def\blem{\begin{lemma}}

\def\elem{\end{lemma}}

\def\bex{\begin{example}}

\def\eex{\end{example}}

\def\bexo{\begin{exercice} \rm }

\def\eexo{\end{exercice} }

\def\brem{\begin{remark}}

\def\erem{\end{remark}}

\def\prf{{\bf Proof }}

\def\bdes{\begin{description}}

\def\edes{\end{description}}

\def\iti{\item[(i)]}

\def\itii{\item[(ii)]}

\def\itiii{\item[(iii)]}

\def\beq{\begin{equation}}

\def\eeq{\end{equation}}

\def\ben{\begin{enumerate}}

\def\een{\end{enumerate}}

\def\beqar{\begin{eqnarray}}

\def\eeqar{\end{eqnarray}}

\def\beqarr{\begin{eqnarray*}}

\def\eeqarr{\end{eqnarray*}}

\def\prf{{\bf Proof }\hspace{.1in}}

\def\Pr{{\mathsf P}}

\def\RR{{\mathbb R}}  

\def\NN{{\mathbb N}}
\def\EE{{\mathbb E}}
\def\PP{{\mathbb P}}

\def\XX{{\mathcal X}}
\def\rar{\rightarrow}
\newcommand{\bi}{\mathbf{i}}
\newcommand{\bu}{\mathbf{u}}
\newcommand{\bv}{\mathbf{v}}
\newcommand{\bphi}{\mathbf{\Phi}}

\def\eps{\varepsilon}
\newcommand{\Bcal}{\mathcal{B}}
\newcommand{\Pcal}{\mathcal{P}}

\def\1{{\rm 1\mskip-4.4mu l}}
\begin{document}
\title{Randomly switched vector fields sharing a  zero on a common invariant face}
\author{Edouard Strickler\\ Institut de Math\'ematiques\\Universit\'e de Neuch\^atel, Switzerland}
\maketitle
\begin{abstract}
We consider a Piecewise Deterministic Markov Process given by random switching between finitely many vector fields vanishing at $0$. It has been shown recently that the behaviour of this process is mainly determined by the signs of Lyapunov exponents. However, results have only been given when all these exponents have the same sign. In this note, we consider the degenerate case where the process leaves invariant some face and results are stated when the Lyapunov exponents are of opposite signs. Applications are given to Lorenz vector fields with switching, and to SIRS model in random environment.

\end{abstract}
\section{Introduction}
In this paper, we consider a Markov process obtained by random switching between finitely many vector fields $F^i : \RR^d \to \RR^d$, sharing a common equilibrium point $q$. It has been shown in \cite{BS17} that the behaviour of the obtained system near $q$ is determined by the sign of quantities linked to classical Lyapunov exponents. These exponents depend on the Jacobian matrices $A^i = DF^i(q)$ and the switching mechanism. There might be $1$ to $d$ distinct exponents. In \cite{BS17}, results are only given in the case where all the Lyapunov exponents have the same sign. Briefly put, it they are all negative, the  system converges to $q$ with positive probability provided the initial condition is close to $q$; while when they are all positive,  the process admits an invariant probability measure that gives no mass to $q$. In the present paper   we consider the degenerate situation, where the process leaves invariant a face  $\{0\} \times \RR^m \subset \RR^d$ containing $q$ so that the Jacobian matrices have the form $$A^i = \begin{pmatrix}
B^i && 0\\
C^i && D^i
\end{pmatrix}.$$ We show that if both maximal Lyapunov exponents associated with $B^i$ and $D^i$ are negative, then the process converges to $q$; while if all the Lyapunov exponents associated to $B^i$ are positive and those to $D^i$ are negative, the process admits an invariant probability measure that gives no mass to $\{0\} \times \RR^m \subset \RR^d$, and hence no mass to $q$. We also notice that in this last case, the Lyapunov exponents associated to $A^i$ take positive and negative values, so that the results of \cite{BS17} cannot be applied. The paper is organised as follows. In section \ref{sec:results}, the main results are stated. The proofs are postponed to section \ref{sec:proofs}. In section \ref{sec:examples},  we give several applications. In particular, our result enables us to close a gap in a discussion on random switching between two Lorenz vector fields in \cite{bakhtin&hurt}. We also recover and slightly extend the results on SIRS models with Markov switching given in \cite{li17}.

\section{Notations and results}
\label{sec:results}
Let $d \geq 1$, $E=\{1,\ldots,N\}$ a finite set and for all $i \in E$, $F^i : \RR^d \to \RR^d$ a $C^2$ globally integrable vector field. We denote by $\varphi^i$ the flow induced by $F^i$ and we assume that there exists a closed set $M$ which is forward invariant for all the vector fields, that is $$\varphi_t^i(M) \subset M, \quad \forall t \geq 0.$$ For all $x \in M$, we are given  an irreducible rate matrix $(a_{ij}(x))_{i,j \in E}$, continuous in $x$. We consider a Markov process $(Z_t)_{t \geq 0} = (X_t, I_t)_{t \geq 0} \in M \times E$, where $X$ evolves according to \beq
 \frac{dX_t}{dt} = F^{I_t}(X_t),
\label{eq:pdmp}
\eeq
and  $I$ is a continuous time jump process taking values in $E$ controlled by $X:$
$$\Pr(I_{t+s} = j | {\cal F}_t, I_t = i)  = a_{ij}(X_t) s + o(s) \mbox{ for } j \neq i \mbox{ on } \{I_t = i\},$$
where ${\cal F}_t = \sigma ((X_s,I_s) \: : s \leq t\}).$
It can be shown (see e.g \cite{BMZIHP}) that the infinitesimal generator of $Z$ is the operator ${\cal L}$ acting on functions $g : M \times E \mapsto \RR,$ smooth in the first variable, according to the formula
\beq
\label{eq:defcL}
{\cal L}g(x,i) =  \langle F^i(x), \nabla g^i(x) \rangle + \sum_{j \in E} a_{ij}(x) (g^j(x) - g^i(x)),
\eeq
where $g^i(x)$ stands for $g(x,i)$. The process $Z$ belongs to the class of  \textit{Piecewise Deterministic Markov Processes} (PDMP), as introduced by Davis in \cite{Dav84}.  

Without loss of generality, we assume that $q=0$. For $n,m$ such that $n+m=d$, and $x \in \RR^d = \RR^n \times \RR^m$, we set $x=(x_n,x_m)$. The notation $0_k$ for $k=n,m$ refers to the zero vector of $\RR^k$. We also write $F^i(x) = (F^i_n(x),F^i_m(x))$. Our standing assumption is :
\begin{hypothesis}\
\label{hyp:stand}
\begin{enumerate}
\item The origin lies in $M$ and for all $i \in E$, $F^i(0)=0$.
\item For all $x_m \in \RR^m$ and all $i \in E$, $F_n^i(0_n,x_m)=0$.
\item The set $M$ intersects the face $\{0_n\} \times \RR^m$ : $\{0\} \subsetneq M \cap ( \{0_n\} \times  \RR^m) \subsetneq M $ 
\item The set $M$  is compact and locally star shaped around the origin : there exists
$\delta > 0$ such that
$$x \in M \mbox{ and } \|x\| \leq \delta \Rightarrow [0,x] \subset M,$$ where $[0,x] = \{ tx, \: t \in [0,1]\}$.
\end{enumerate}
\end{hypothesis}
The second assumption implies that the face $\{0_n\} \times \RR^m$ is invariant under each $\varphi^i$ : for all $t \geq 0$, $\varphi_t^i(x) \in \{0_n\} \times  \RR^m $ if and only if $x \in \{0_n\} \times  \RR^m$. We set $M_+ = \{ (x_n,x_m) \in M \: : \: x_n \neq 0 \}$ and $M_0 = M \setminus M_+$. Both $M_0$ and $M_+$ are non empty, and $M_0$ is invariant for all the flows $\varphi^i$.  For all $i \in E$, set $A^i = DF^i(0)$, the Jacobian matrix of $F^i$ at $0$. The second assumption has also the consequence that $A^i$ is block lower triangular : 
\beq
\label{eq:triang}
A^i = \begin{pmatrix}
B^i && 0\\
C^i && D^i
\end{pmatrix},
\eeq
with $B^i \in M_n(\RR)$, $C^i \in M_{m,n}(\RR)$ and $D^i \in M_m(\RR)$. 

\subsection{Notation}
Throughout the paper we will adopt the following notation : $\langle \cdot, \cdot  \rangle$ denotes the Euclidean inner product in $\RR^k,$ for $k=n,m,d$;  $\| \cdot \|$ the associated norm, and $S^{k-1} = \{x \in \RR^k \: : \|x\| = 1\}$ the unit sphere. 
 For a metric space $(\XX,d)$, we will denote by $\Bcal(\XX)$ the set of Borel sets of $\XX$, and by $\Pcal(\XX)$ the set of probability measures on $\Bcal(\XX)$. If $(Z_t)_{t \geq 0}$ is a Markov process on $\XX$ and $\nu \in \Pcal(\XX)$, we set, as usual, $\PP^Z_{\nu}$ for the law of the process $Z$ with initial distribution $\nu$ and $\EE^Z_{\nu}$ for the associated expectation. If $\nu = \delta_x$ for some $x \in \XX$, we write $\PP^Z_x$ for $\PP^Z_{\delta_x}$.  We denote by $(P^Z_t)_{t \geq 0}$ the semigroup of $Z$ acting on bounded measurable function $f : \XX \to \RR$ as $$ P^Z_t f(x) = \EE_x^Z \left( f(Z_t) \right).$$ When there is no ambiguity on the process considered, we drop the exponent $Z$. An invariant distribution for the process $Z$ is a probability $\mu \in \Pcal(\XX)$ such that $\mu P_t = \mu$ for all $t \geq 0$. We let $\Pcal_{inv}^Z$ denote the set of all the invariant distributions of $Z$ and for $N \subset \XX$, let  $\Pcal_{inv}^Z(N)$ denote the (possibly empty) set of invariant probabilities giving mass 1 to the set $N$.  For  $\bi=(i_1, \ldots, i_k) \in E^k$ and $\bu = (u_1,\ldots,u_k) \in \RR_+^k$, we denote by $\bphi_{\bu}^{\bi}$ the composite flow : $\bphi_{\bu}^{\bi} = \varphi_{u_k}^{i_k} \circ \ldots \circ \varphi_{u_1}^{i_1}$.  For $x \in M$ and $t \geq 0$, we denote by $\gamma^+_t(x)$ (resp. $\gamma^+(x)$)  the set of points that are reachable from $x$ at time $t$ (resp. at any nonnegative time) with a composite flow: $$\gamma^+_t(x)=\{ \bphi_{\bv}^{\bi}(x), \: (\bi,\bv) \in E^k \times \RR_+^k, k \in \NN, v_1 + \ldots + v_k = t\},$$

$$ 
\gamma^+(x) = \bigcup_{t \geq 0} \gamma^+_t(x).
$$ 
We will say that a point $x^* \in M$ is \emph{accessible} from $B \subset M$ if $x^* \in \cap_{x \in B} \overline{\gamma^+(x)}$.

\subsection{Linear system and Lyapunov exponents}
\label{sec:linear}

For a given set of matrices $\hat{A}=(\hat{A}^i)_{i \in E}$ of size $k \times k$, we consider the linear system $(Y,J)$ where $Y$ evolves according to $$\frac{d Y_t}{d t} = \hat{A}^{J_t} Y_t,$$ and $J$ is a continuous time Markov chain on $E$ with transition rate matrix $(a_{ij}(0))_{i,j \in E}$. By irreducibility of $(a_{ij}(0))_{i,j \in E}$, $J$ admits a unique invariant probability measure on $E$ denoted by $p$.


Whenever the initial condition $y_0$ is not zero, the angular part of $Y_t$, $\Theta_t = \frac{Y_t}{\|Y_t\|}$ is well defined, and evolves according to 

\beq
\label{dThetarho}
\frac{d\Theta_t}{dt} = \hat{A}^{J_t}\Theta_t -  \langle A^{J_t}\Theta_t, \Theta_t \rangle \Theta_t.
 \eeq
This defines a differential equation on 
 $S^{k-1}$ and the process $(\Theta_t,J_t)_{t \geq 0}$ is a PDMP on $S^{k-1} \times E$. When we need to emphasis the dependence on $(\hat{A}^i)_{i \in E}$, we denote by $\Theta(\hat{A})$ the solution of \eqref{dThetarho}. For an invariant probability $\mu$ of $(\Theta(\hat{A}),J)$, we define the $\mu$-{\em average growth rate} as
\beq
\label{deflambdamu}
\Lambda_{\hat{A}}(\mu) = \int \langle \hat{A}^i \theta, \theta \rangle \mu(d\theta di) = \sum_{i \in E} \int_{S^{k-1}} \langle \hat{A}^i \theta, \theta \rangle \mu^i(d\theta),
\eeq where
$\mu^i(\cdot)$ is the measure on $S^{k-1}$ defined by $$\mu^i(\cdot) = \mu (\cdot \times \{i\}).$$
We let $\Lambda(\hat{A})$ be the set of all the $\Lambda_{\hat{A}}(\mu)$ for $\mu$ invariant probability of $(\Theta(\hat{A}),J)$. As in \cite{BS17}, we define the {\em extremal average growth rates} as the numbers
\beq
\label{deflambdamumax}
\Lambda_{\hat{A}}^{-} = \inf \Lambda(\hat{A}) \mbox{ and }
\Lambda_{\hat{A}}^+ = \sup \Lambda(\hat{A}).
\eeq
In \cite{BS17}, we show that  $\Lambda(\hat{A})$ is composed of Lyapunov exponents in the sense  of Oseledet's Multiplicative Ergodic Theorem (see e.g \cite[Theorem 3.4.1]{Arn98} and Section \ref{sec:proofs}). In particular, $\Lambda(\hat{A})$ is actually a finite set, and  the supremum and the infimum in equation \eqref{deflambdamumax} are maximum and minimum. 
We start with a lemma.

\blem
\label{lem:triang}
Assume that all the $A^i$ have the block triangular form \eqref{eq:triang}.Then, with the above notations, $\Lambda(D) \subset \Lambda(A)$.  
\elem

\prf
Let $\lambda \in \Lambda(D)$ and $\hat{\mu}$ be an invariant probability of $(\Theta(D),J)$ such that $\lambda = \Lambda_D(\hat{\mu})$. For $\Theta \in S^{d-1}$, we write $\Theta=(\Theta^n,\Theta^m)$. We this notation, \eqref{dThetarho} becomes :

\beq
\label{dThetarnm}
 \left \{
\begin{array}{l}
 \frac{d \Theta^n_t}{dt} = B^{J_t} \Theta^n_t -\left( \langle B^{J_t} \Theta^n_t,  \Theta^n_t \rangle +  \langle C^{J_t} \Theta^n_t + D^{J_t} \Theta^m_t, \Theta^m_t\rangle\right) \Theta^n_t \\ \\
 \frac{d \Theta^m_t}{dt} =C^{J_t} \Theta^n_t + D^{J_t} \Theta^m_t -\left( \langle B^{J_t} \Theta^n_t,  \Theta^n_t \rangle +  \langle C^{J_t} \Theta^n_t + D^{J_t} \Theta^m_t , \Theta^m_t \rangle \right) \Theta^m_t
\end{array}
\right.
\eeq
From this equation, one can see that the space $\{(\theta_n,\theta_m) \in S^{d-1} \: : \: \theta_n=0\}$ is invariant, and on that space,$(\Theta(A),J)=(0,\Theta(D),J)$. Now  we extend $\hat{\mu}$ to a probability measure $\mu$ on $S^{d-1}\times E$ such that $\mu(\{(\theta_n,\theta_m) \in S^{d-1} \: : \: \theta_n=0\} \times E)=1$ and the marginal of $\mu$ on $S^{m-1} \times E$ is $\hat{\mu}$. Then, $\mu$ is an invariant probability for $(\Theta(A),J)$, and straightforward computation shows that $\Lambda_A(\mu)=\Lambda_D(\hat{\mu})=\lambda$. Thus $\lambda \in \Lambda(A)$.
\QEDB
\bigskip

\brem {\rm
The same proof shows that in case where the $A^i$ are block diagonal, that is $C^i = 0$, then $\Lambda(B) \subset \Lambda(A)$. However, this is not true in general. Here is a counter  example in dimension $d=2$. Let $A^i$, $i=0,1$ be two $2 \times 2$ matrices defined by $$ A^i= \begin{pmatrix}
b_i & 0 \\
c_i & d_i
\end{pmatrix},$$
and assume that $b_i < d_i$ for $i=0,1$ as well as $c_0(b_1-d_1) \neq c_1 (b_0 - d_0)$. In particular, $\Lambda_B^+=\sum_ip_ib_i < \sum_i p_i d_i = \Lambda^-_D$.  We show that in this case, the set of invariant probability measures of $(\Theta(A),J)$ reduces to $\delta_{(0,1)} \otimes p$ and $\delta_{(0,-1)} \otimes p$; hence $\Lambda(A) = \Lambda(D) \nsupseteq \Lambda(B)$.   Let $\theta_i$ be the normalized eigenvector of $A^i$ associated with $b_i$. Since  $c_0(b_1-d_1) \neq c_1 (b_0 - d_0)$, $\theta_0 \neq \theta_1$. Now it is easily checked that the region between $\theta_0$ and $\theta_1$ is transient for $\Theta(A)$ and that when the process leaves this region, $\Theta_t(A)$ converges to $(0,1)$ or $(0,-1)$. 
} \erem 
 
We prove in Section \ref{sec:proofs} that the result given in the preceding remark can be generalized as follows :

\bprop
\label{prop:lambdaA}
Assume that all the $A^i$ have the block triangular form \eqref{eq:triang}. If $\Lambda_B^+ < \Lambda_D^-$ and if $\{(\theta_n,\theta_m) \in S^{d-1} \: : \: \theta_n=0\}$ is accessible from $S^{d-1}$, then $\Lambda(A) = \Lambda(D)$.
\eprop
Using a result from Hennion \cite{H84}, we have the following proposition, whose proof is given in Section \ref{sec:proofs}.

\bprop
\label{prop:triang}
Assume that all the $A^i$ have the block triangular form \eqref{eq:triang}.Then, with the above notations, $\Lambda_A^+ = \max( \Lambda_B^+, \Lambda_D^+).$
\eprop

\bex
\label{ex:dim2}
We describe completely the two dimensional case. Let $(A^i)_{i \in E}$ be a family of $2 \times 2$  upper triangular matrices : 

$$ A^i= \begin{pmatrix}
b_i & 0 \\
c_i & d_i
\end{pmatrix}.$$

One has $\Lambda_B^+=\Lambda_B^-= \sum_i p_i b_i:=\Lambda_B $ and $\Lambda_D^+=\Lambda_D^-= \sum_i p_i d_i:=\Lambda_D$. We have the following :
\begin{enumerate}
\item If $\Lambda_B > \Lambda_D$, then $\Lambda_A^+ = \Lambda_B$ and $\Lambda_A^-=\Lambda_D$;
\item If $\Lambda_B = \Lambda_D$, then $\Lambda_A^+ = \Lambda_A^-= \Lambda_B=\Lambda_D $;
\item If for all $i \neq j$, $c_i (b_j-d_j)=c_j (b_i - d_i)$, then $\Lambda_A^+ = \max(\Lambda_B,\Lambda_D)$ and $\Lambda_A^- = \min(\Lambda_B,\Lambda_D)$;
\item If $\Lambda_B < \Lambda_D$ and if there exist $i \neq j$ such that $c_i (b_j-d_j) \neq c_j (b_i - d_i)$, then $\Lambda_A^+ = \Lambda_A^-= \Lambda_D$.
\end{enumerate}

In case where $\Lambda_B > \Lambda_D$, then $\Lambda_A^+= \Lambda_B$ by Proposition \ref{prop:triang} and since by Lemma \ref{lem:triang}, $\Lambda(D) \subset \Lambda(A)$, one has $\Lambda_A^-=\Lambda_D$. 

If $\Lambda_B = \Lambda_D$, the set of Lyapunov exponents of $(A^i)_{i \in E}$ in the sense of ergodic theory reduces to $\Lambda_B$, hence the result (see proof of Proposition \ref{prop:lambdaA} in Section \ref{sec:proofs}).

Now assume that for all $i \neq j$, $c_i (b_j-d_j)=c_j (b_i - d_i)$. If $\Lambda_B = \Lambda_D$, then the result follows from point 2. If  $\Lambda_B \neq \Lambda_D$, there exists $i_0 \in E$ such that $b_{i_0} \neq d_{i_0}$. Set $x^*=(1, \frac{c_{i_0}}{b_{i_0} - d_{i_0}})$  We claim that $x^*$ is a common eigenvector for all the $A^i$. Indeed, let $i \in E$. If $b_i-d_i \neq 0$, then $x=(1,\frac{c_i}{b_i - d_i})$ is an eigenvector of $A^i$ associated with $b_i$, and since $c_i (b_{i_0}-d_{i_0})=c_{i_0} (b_i - d_i)$, $x=x^*$. If $b_i-d_i = 0$, since $b_{i_0}-d_{i_0} \neq 0$ and $c_i (b_{i_0}-d_{i_0})=c_{i_0} (b_i - d_i)$, one has $c_i=0$, in other words $A^i = b_i I$, hence $x^*$ is an eigenvector of $A^i$. We conclude that if we let $\theta^* = \frac{x^*}{\|x^*\|},$ then $\mu=\delta_{\theta^*} \otimes p \in \Pcal_{inv}^{(\Theta(A),J)}$ and $\Lambda_A(\mu) = \Lambda_B$. This combined with Lemma \ref{lem:triang} proves point 3.

Finally assume that $\Lambda_B < \Lambda_D$. It implies that there exists $i_0 \in E$ such that $b_{i_0} < d_{i_0}$. Thus, $\{(0,1),(0-1)\}$ is accessible from every point in $S^{d-1} \setminus \{ \theta^* \}$ where $\theta^*$ is defined as before. Now there exist $j \in E$ such that $\theta^*$ is not an eigenvector for $A^j$. In particular, we can reach $S^{d-1} \setminus \{ \theta^* \}$ from $\theta^*$ by following $A^j$. Hence  $\{(0,1),(0-1)\}$ is accessible from $S^{d-1}$ and the result follows from Proposition \ref{prop:lambdaA}.

\eex

\subsection{Main Results}
The first theorem is an immediate consequence of Proposition \ref{prop:triang} and Theorem 3.1 in \cite{BS17}. 

\bthm
\label{thm:ext}
Assume $\Lambda_B^+<0$ and $\Lambda_D^+<0$. Let $0 < \alpha < - \Lambda_A^+.$ Then there exists a neighborhood ${\cal U}$ of $0$  and $\eta > 0$ such that for all $x \in {\cal U}$ and $i \in E$
$$\mathbb{P}_{x,i}( \limsup_{t \rar \infty} \frac{1}{t} \log(\|X_t\|) \leq - \alpha) \geq \eta.$$
If furthermore $0$ is accessible from $M,$ then for all $x \in M$ and $i \in E$
$$\mathbb{P}_{x,i}( \limsup_{t \rar \infty} \frac{1}{t} \log(\|X_t\|) \leq \Lambda_A^+) = 1.$$
\ethm

The next theorem is the main result of this paper. It gives results when the Lyapunov exponents are of opposite signs. We let $$\Pi_t = \frac{1}{t} \int_0^t \delta_{Z_s} ds \in {\cal P}(M \times E)$$ denote  the {\em empirical occupation measure} of the process $Z.$ For every Borel set $A \subset M \times E$
$$\Pi_t(A) = \frac{1}{t} \int_0^t \1_{\{Z_s \in A\}} ds$$ is then the proportion of the time spent by $Z$ in $A$ up to time $t.$ Recall that $M_+=  \{ (x_n,x_m) \in M \: : \: x_n \neq 0 \}$, and for $\delta > 0$, set $M_0^{\delta} = \{ (x_n,x_m) \in M_+ \: : \: \|x_n\| < \delta \}.$

\bthm
\label{thm:pers}
Assume $\Lambda_B^- > 0 > \Lambda_D^+$
 and that $0$ is accessible from $M_0 \times E$. Then :
\bdes \iti For all $\eps > 0$ there exists $\delta > 0$ such that for all $x \in M_+$, $i \in E$, $\mathbb{P}_{x,i}$ almost surely,
$$\limsup_{t \rar \infty} \Pi_t(M_0^{\delta} \times E) \leq \eps.$$
In particular, for all $x \in M_+,$ $\mathbb{P}_{x,i}$ almost surely,  every limit point (for the weak* topology) of $(\Pi_t)$ belongs to ${\cal P}_{inv} ( M_+ \times E).$
\itii There exist positive constants $\theta, K$ such that for all $\mu \in {\cal P}_{inv}( M_+ \times E)$
$$\sum_{i \in E} \int \|x_n\|^{-\theta} d \mu^i(x_n,x_m) \leq K.$$
\itiii Let $\eps > 0$ and $\tau^{\eps}$ be the stopping time defined by
$$\tau^{\eps} = \inf \{t \geq 0 : \: \|X_t^n\| \geq \eps \}.$$ There exist $\eps > 0$,  $b > 1$ and $c > 0$ such that
for all $x \in M_+$ and $i \in E$,
$$\mathbb{E}_{x,i}^Z(b^{\tau^{\eps}}) \leq c (1 + \|x_n\|^{-\theta}).$$
\edes
\ethm

\brem
{\rm Note that under the assumptions of the above theorem,  Lemma \ref{lem:triang} implies $\Lambda^-_A \leq \Lambda^-_D < 0$ while by Proposition \ref{prop:triang}, $\Lambda_A^+ = \Lambda_B^+ > 0$. Thus the results in \cite{BS17} cannot be applied.

}\erem
As in \cite{BS17}, we give the following theorem ensuring uniqueness of the invariant probability giving no mass to $M_0 \times E$ which is a consequence of results in \cite{BMZIHP} (see also \cite{bakhtin&hurt}). Set $\mathrm{F}_0 = \{F^i\}_{i \in E}$ and $\mathrm{F}_{k+1} = \mathrm{F}_k \cup \{[F^i, V], V \in \mathrm{F}_k\}$ where $[, ]$ is the Lie bracket operation.
 We say  that  the {\em weak bracket} condition holds at $x \in M$ provided the vector space spanned by the vectors $\{V(x) \: : V \in \cup_{k \geq 0}  \mathrm{F}_{k}\}$  has full rank.
 We let $\mathbf{Leb}$ denote the Lebesgue measure on $\RR^d.$
\bthm
 \label{th:persist2} In addition to the assumptions of Theorem \ref{thm:pers},  suppose that there exists a point $y \in M_+$ accessible from $M_+$ at which the weak bracket condition holds. Then
 \bdes
 \iti
The set ${\cal P}_{inv}^Z ( M_+ \times E)$ reduces to a single element, denoted $\Pi$;
\itii  $\Pi$ is absolutely continuous with respect to $\mathbf{Leb} \otimes (\sum_{i \in E} \delta_i)$;
\itiii For all $x \in M_+$ and $i \in E$, $$\lim_{t \rar \infty} \Pi_t = \Pi$$ $\mathbb{P}^Z_{x,i}$ almost surely. \edes
\ethm

Set  $\mathcal{F}_{0} = \{F^i  - F^j\: : i, j = 1, \ldots N\}$ and $\mathcal{F}_{k+1}  =  \mathcal{F}_k \cup \{[F^i , V ] \: : V \in \mathcal{F}_k\}.$
We say that the {\em strong bracket} condition holds at $x \in M$ provided the vector space spanned by the vectors $\{V(x) \: : V \in \cup_{k \geq 0}  \mathcal{F}_{k}\}$  has full rank.

\bthm
 \label{th:persist3}
 In addition to the assumptions of Theorem \ref{thm:pers},  suppose that one of the two following holds :
\bdes
\iti The weak bracket condition is strengthened to the strong bracket condition; or
\itii There exist $\alpha_1,\ldots,\alpha_N \in \RR$ with $\sum \alpha_i = 1$ and a point $e^{\star} \in M_+$ accessible from $M_+$ such that  $\sum \alpha_i F^i(e^{\star}) = 0$.
\edes
  Then there exist $\kappa, \theta > 0$ such that for all $x \in M_+$ and $i \in E$,
$$\|\mathbb{P}_{x,i}(Z_t \in \cdot ) - \Pi\|_{TV} = \| \delta_{x,i} P_t^Z - \Pi \|_{TV} \leq const. (1 + \|x_n\|^{-\theta}) e^{-\kappa t}.$$
\ethm

\section{Examples}
\label{sec:examples}
In this section we give several examples of applications of our results.
\subsection{Lorenz Vector Fields}
In \cite{bakhtin&hurt}, the authors consider a random switching between two Lorenz vector fields $F^i$, $i=0,1$ :
\beq
\label{eq:lorenz}
F^i(x,y,z) = \begin{pmatrix}
\sigma_i ( y -x)\\
r_i x - y - xz\\
xy - b_i z
\end{pmatrix},
\eeq
with $\sigma_0=\sigma_1=10$, $b_0=b_1=8/3$, $r_0=28$, and $r_1 \neq r_0$ close to $28$. It is known since the proof of Tucker \cite{T99} that $F^0$ admits a robust strange attractor $\Gamma_0$. Thus for $r_1$ close to $r_0$, $F^1$ shares this property.  
In \cite{bakhtin&hurt}, it has been shown that the compact set   $M = \{(x,y,z) \in \RR^3 \: : \: 2r_0 \sigma (x^2+y^2)+2 \sigma b (z_0 - r_0)^2 \leq 2 \sigma b r_0^2\}$  is forward invariant, and that $\Gamma_0$ is accessible from every point that does not lie on the $z$-axis. Moreover they proved that the strong bracket condition holds at any point which is not on the $z$-axis. Then they argue that by compactness of $M$, there exists an invariant probability, and that it has to be absolutely continuous with respect to the Lebesgue measure due to the bracket condition. However, this argument is not sufficient : there exists indeed an invariant probability measure on $M$, which is $\delta_0 \otimes p$. However, this measure is not absolutely continuous. We explain how our results apply to that situation and fill in this gap in the proof of \cite{bakhtin&hurt}. In particular, we prove the following result :

\bprop
Let $F^i$, $i=0,1$ be two Lorenz vector fields defined by \eqref{eq:lorenz} with $\sigma_0=\sigma_1=10$, $b_0=b_1=8/3$, $r_0=28$, and $r_1 \neq r_0$ close to $28$. Then the process $Z$ admits a unique invariant probability measure $\Pi$ such that $\Pi(M \setminus \{ x = y = 0\}) = 1$. Moreover, $\Pi$ is absolutely continuous with respect to the Lebesgue measure, and there exist $\kappa, \theta > 0$ such that for all $x_0=(x,y,z) \in M$ such that $(x,y) \neq 0$ and $i \in E$,
$$\|\mathbb{P}^Z_{x,i}(Z_t \in \cdot ) - \Pi\|_{TV}  \leq const. (1 + \|(x,y)\|^{-\theta}) e^{-\kappa t}.$$
\eprop 

\prf
One can note that the $z$-axis is invariant, and that assumption \ref{hyp:stand} holds with $n=2$ and $m=1$. Moreover, we have 
\beq
A^i = \begin{pmatrix}
B^i & 0\\
0 & - b_i
\end{pmatrix},
\eeq
where $$B^i = \begin{pmatrix}
- \sigma_i & \sigma_i \\
r_i & -1
\end{pmatrix}.$$
Setting $D^i=(-b_i)$, one has $\Lambda^+_D = - ( p_0 b_0 + p_1 b_1) < 0$. Furthermore, on the $z$ - axis, $|z_t| \leq z_0 e^{- b t}$, with $b = \min(b_0,b_1)$. Hence $0$ is accessible from $M_0$. Let us show that $\Lambda^-_B > 0$. First, it is easily checked that $B^0$ and $B^1$ have no common eigenvectors. Therefore, \cite[Example 2.12]{BS17} implies that $\Lambda_B^+=\Lambda_B^-$. Next, the matrices $B^i$ are Metzler, meaning that their off diagonal entries are nonnegative. Therefore, the Kolotilina-type lower estimate for the
top Lyapunov exponent proved by Mierczy{\'n}ski \cite[Theorem 1.3]{Mier13} implies that $$ \Lambda_B^-  \geq \frac{1}{2} \sum_i p_i Tr(B^i) +  \sum_i p_i \sqrt{B_{12}^i B_{21}^i}.$$
Here, $Tr(B^i) = - 11$ for $i=0,1$, and $\sqrt{B_{12}^0 B_{21}^0} = \sqrt{280} > 11/2$. Since $r_1$ is close to $r_0$, we also have that $\sqrt{B_{12}^1 B_{21}^1} > 11/2$, hence $\Lambda_B^- > 0$. The result follows from Theorem \ref{th:persist3} due to the strong bracket condition proved in \cite{bakhtin&hurt}.  \QEDB   
\bigskip

\begin{figure}[!h]
\centering
\includegraphics[scale=0.3]{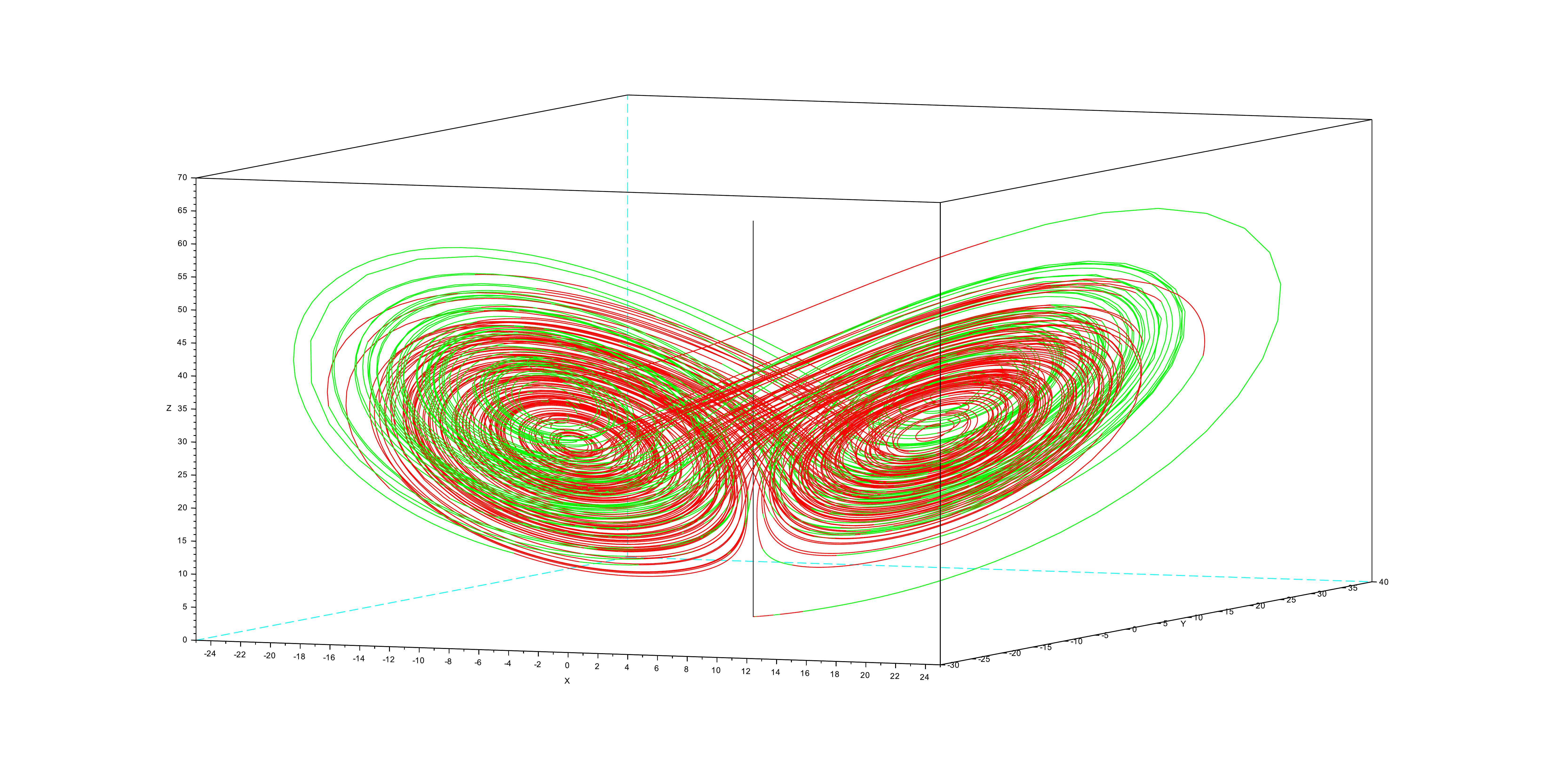}
\caption{\label{fig:lorenz} Randomly switched Lorenz vector fields}
\end{figure}

In Figure \ref{fig:lorenz}, we show a trajectory of $X_t$ with initial condition $(0,0.05,0.05)$ for the vector fields $F^0$ and $F^1$ given by the above values of parameters and $r_1=35$. The $z$-axis is drawn in black. 

\subsection{Epidemiological SIRS models}
In this section, we show  how our result enables to recover and extend those found in \cite{li17}. In this paper, the following SIRS model with random switching is studied  :

\begin{equation}\label{e:SIRS}
F^k(S,I,R) =  \begin{pmatrix}
\Lambda- \mu S + \lambda_k R - \beta_k S G_k(I)\\
\beta_k S G_k(I)-(\mu+\alpha_k+\delta_k)I\\
\delta_k I -(\mu + \lambda_k)R
\end{pmatrix},
\end{equation} 
for $k \in E=\{1,\ldots,N\}$, where $G_k$ is a regular function such that $G_k(0)=0$.  The reader is referred to \cite{li17} for the epidemiological interpretation of the different constants. The authors study the specific case where only $\beta$ is allowed to depend on $k$ and where the discrete component $(r_t)_{t \geq 0}$ is an irreducible Markov chain on $E$, that is the rate matrix $a$ does not depend on the position. Here we assume that the positive constants $\lambda_k, \alpha_k, \delta_k$ and the functions $G_k$ may depends on $k$ and that $a$ could depend on the position. We still let $\Lambda$ and $\mu$ be constant : they are the intrinsic birth and death rates, and are not related on how the disease spread among the population. Thus the point $q=(\frac{\Lambda}{\mu},0,0)$ is a common equilibrium for the $F^k$.  Set $(Z_t)_{ t \geq 0} = (X_t,r_t)_{ t \geq 0}$, with $X_t=(S_t,I_t,R_t)$.  Write $\RR_{+}^3 = \{ x \in \RR^3 \: : x_i \geq 0, i=1,2,3\}$ and $\RR_{++}^3 = \{ x \in \RR^3 \: : x_i > 0, i=1,2,3\}$. It is easily seen that $\RR_+^3$ and $\RR_{++}^3$ are positively invariant for all the $F^k$. Moreover, one can check that the compact set $$M = \{ x \in \RR_+^3 \: : \: x_1+x_2+x_3 \leq \Lambda/\mu \}$$ is also positively invariant for all the $F^k$. Furthermore, there are two invariant sets : the $S$-axis and the $(S,R)$ - plane. We set $M_0 = \{ (S,I,R) \in M \: : \: I= 0\}$.
  We make the following assumptions, that are taken from \cite{li17} :

\begin{hypothesis}\
\label{hyp:SIR}
\bdes
\iti For all $k$, $G_k : \RR_+ \to \RR_+$ is $C^2$, with $G_k(0)=0$ and $0 < G_k(I) \leq G_k'(0)I$ for $I > 0$;
\itii For all $k$, if $\beta_k \frac{\Lambda}{\mu} G_k'(0) - (\mu + \alpha_k + \delta_k) > 0$, then $F^k$ admits an equilibrium point $x^* \in M_+$ which is accessible from $M_+$.
\edes
\end{hypothesis}

 For convenience, we reorder the coordinates as $(I,R,S)$ and set $q=(0,0,\frac{\Lambda}{\mu})$. Writing $A^k=DF^k(q)$, one has 

$$ A^k = \begin{pmatrix}
\beta_k \frac{\Lambda}{\mu} G_k'(0) - (\mu + \alpha_k + \delta_k) & 0 & 0\\
\delta & -(\mu+\lambda_k) & 0\\
- \beta_k \frac{\Lambda}{\mu} G_k'(0) & \lambda_k & - \mu
\end{pmatrix}.$$
If we denote by $D^k$ the matrix $$ D^k = \begin{pmatrix}  -(\mu+\lambda_k) & 0\\
 \lambda_k & - \mu
\end{pmatrix},$$
 then by Proposition \ref{prop:triang}, $\Lambda_D^+= \max( \Lambda_1, \Lambda_2)$, with
  
  $$\Lambda_1 = - \sum_k p_k (\mu + \lambda_k) < 0$$ 
  and 
  
  $$\Lambda_2 = -  \sum_k  p_k \mu =  - \mu <0.$$
   Hence $\Lambda_D^+ = - \mu < 0$, and by Theorem \ref{thm:ext}, on $M_0$, the process converges to $q$. Now if $B^k = (\beta_k \frac{\Lambda}{\mu} G'(0)-(\mu + \alpha_k + \delta_k))$, then $\Lambda_B^- = \Lambda_B^+ = \sum_k p_k (\beta_k \frac{\Lambda}{\mu} G_k'(0)-(\mu + \alpha_k + \delta_k))$. As in \cite{li17}, we set  $$ R_0 = \frac{\sum_k p_k \beta_k \frac{\Lambda}{\mu} G'(0)}{\sum_k p_k (\mu + \alpha_k + \delta_k)}.$$ Note that $R_0 < 1$ (respectively $R_0>1$) if and only if $\Lambda_B^- < 0$ (resp. $\Lambda^-_B>0$). In particular, Theorems   \ref{thm:ext}  and  \ref{thm:pers} imply the following statement, that recovers and slightly extends Theorems 4, 8 and 9 in \cite{li17}.

\bthm
With the above notation, the following hold.
\bdes
\iti Assume that $R_0 < 1$. Then, for all $z_0=(s_0,i_0,r_0,k_0) \in M \times E$, one has $$\PP_{z_0}(\limsup_{t \rar \infty} \frac{1}{t} \log(\|(S_t,I_t,R_t)-(\frac{\Lambda}{\mu},0,0)\|) \leq \Lambda_A^+) = 1,$$ where $\Lambda_A^+ = \max(\Lambda_B^+,-\mu)$.
\itii Assume that $R_0 > 1$. Then the process $Z$ admits an invariant probability measure $\Pi$ such that $\Pi(M \setminus M_0^1 \times E) = 1$. 
\itiii Assume in addition to $R_0>1$ that the weak bracket condition holds at an accessible point. Then $\Pi$ is unique and there exist $\kappa, \theta > 0$ such that for all $x=(s,i,r) \in M_+$ and $k \in E$,
$$\|\mathbb{P}_{x,k}(Z_t \in \cdot ) - \Pi\|_{TV} \leq const. (1 + \|i\|^{-\theta}) e^{-\kappa t}.$$ In addition,  for all $x \in M_+$ and $k \in E$, $$\lim_{t \rar \infty} \Pi_t = \Pi$$ $\mathbb{P}_{x,k}$ almost surely.  
\edes
\ethm

\prf
If $R_0 < 1$, then $\Lambda_B^+<0$ and thus there exists $k_0 \in E$ such that  $\beta_{k_0} \frac{\Lambda}{\mu} G_{k_0}'(0)-(\mu + \alpha_{k_0} + \delta_{k_0}) < 0$. We show that this implies that $q$ is accessible from $M_+$. Let $x_0 \in M_+$ and denote by $x_t=(s_t,i_t,r_t)$ the solution of $$ \frac{d x_t}{dt}=F^{k_0}(x_t)$$ with initial condition $x_0$. Now by assumption \ref{hyp:SIR} and the fact that $s_t \leq \frac{\Lambda}{\mu}$,
 $$\frac{d i_t}{d t} \leq \left( \beta_{k_0} \frac{\Lambda}{\mu} G_{k_0}'(0)-(\mu + \alpha_{k_0} + \delta_{k_0}) \right) i_t.$$ Since $\beta_{k_0} \frac{\Lambda}{\mu} G_{k_0}'(0)-(\mu + \alpha_{k_0} + \delta_{k_0}) < 0$, $i_t$ converges to $0$ exponentially fast. It is easy to check that on $M_0$, $(s_t,r_t)$ converges to $(\frac{\Lambda}{\mu},0)$, thus $x_t$ converges to $q$. Hence $q$ is accessible, and \textbf{(i)} follows from Theorem \ref{thm:ext}. Point \textbf{(ii)} is an immediate consequence of Theorem \ref{thm:pers}. Now if $R_0 > 1$, there exists $k_0 \in E$ such that  $\beta_{k_0} \frac{\Lambda}{\mu} G_{k_0}'(0)-(\mu + \alpha_{k_0} + \delta_{k_0}) > 0$. By assumption \ref{hyp:SIR}, this implies that $F^{k_0}$ admits an accessible equilibrium $x^* \in M_+$. Point \textbf{(iii)} follows then by Theorems \ref{th:persist2} and \ref{th:persist3}. \QEDB

\section{Proofs}
\label{sec:proofs}
\subsection{Proof of Theorem \ref{thm:pers}}
The idea of the proof is similar to that used in \cite{BS17}, and also relies on results of \cite{B18}. In \cite{BS17}, we rewrite the process in spheric coordinates on $\RR_+ \times S^{d-1}$. Here the idea is to only write the spheric coordinates for the part of $X_t$ living in $\RR^n$. That is, we consider the map $ \Psi : \RR^n \setminus \{0_n \} \times \RR^m \times E \to \RR_+^* \times S^{n-1} \times \RR^m \times E$ defined by $ \Psi(x_n,x_m,i) = ( \|x_n\|, \frac{x_n}{\|x_n\|}, x_m,i)$. We set $\mathcal{X}_+ = \Psi(M_+ \times E)$. When $(x,i) \in M_+ \times E$, the process $\tilde{Z}_t = \Psi(Z_t) = (\rho_t,\Theta_t,X^m_t,I_t)$ is well defined  and satisfies 
\beq
\label{dThetarhoNL}
 \left \{
\begin{array}{l}
 \frac{d \rho_t}{dt} =  \langle \Theta_t, \tilde{F}_n^{I_t}(\rho_t,\Theta_t,X^m_t) \rangle \rho_t \\ \\
 \frac{d\Theta_t}{dt} = \tilde{F}_n^{I_t}(\rho_t,\Theta_t,X^m_t)-\langle \Theta_t, \tilde{F}_n^{I_t}(\rho_t,\Theta_t,X_t^m) \rangle \Theta_t \\ \\
  \frac{d X^m_t}{dt} = \tilde{F}_m^{I_t}(\rho_t,\Theta_t,X^m_t) \\ \\
  \Pr(I_{t+s} = j | {\cal F}_t) = a_{i j}(\rho_t \Theta_t, X^m_t) s + o(s) \mbox{ for } i \neq j \mbox{ on } \{I_t = i \}
\end{array}
\right.
\eeq
where for all $(\rho, \theta, x_m) \in \RR_+^* \times S^{n-1} \times \RR^m$, $\tilde{F}_m^i( \rho, \theta, x_m) = F^i_m( \rho \theta, x_m)$ and 
$$\tilde{F}_n^i( \rho, \theta, x_m) = \frac{F^i_n( \rho \theta, x_m)}{\rho}.$$ Since $F^i_n$ is $C^2$ and $F^i_n(0,x_m)=0$, the map $\tilde{F}_n^i$ extends to a $C^1$ map on $\RR_+ \times S^{n-1} \times \RR^m$ by setting $$ \tilde{F}_n^i( 0, \theta, x_m) = B^i(x_m) \theta,$$ where $B^i(x_m) \in M_n(\RR)$ is such that $DF^i_n(0,x_m) = ( B^i(x_m), 0)$. Note that in particular, $B^i(0_m) = B^i$. Thanks to this definition, we can extend \eqref{dThetarhoNL} to $$\mathcal{X} := \overline{\mathcal{X}_+} = \mathcal{X}_+ \cup \mathcal{X}_0$$
where $\mathcal{X}_0 = \{0\} \times S^{n-1} \times \RR^m \times E$. This induces a PDMP  (still denoted $\tilde{Z}$) on $\mathcal{X},$ whose infinitesimal generator $\tilde{\mathcal{L}}$ acts on functions $f : \mathcal{X} \to \RR$ smooth in $(\rho, \theta, x_m)$ according to

\begin{equation}
\label{def:Lpolaire}
\begin{split}
      \tilde{\mathcal{L}} f(\rho, \theta, x_m, i) =& \frac{\partial f^i}{\partial \rho}(\rho, \theta, x_m) \langle \theta, \tilde{F}_n^i(\rho,\theta,x_m) \rangle \rho + \langle \nabla_{\theta} f^i (\rho, \theta, x_m), \tilde{G}^i(\rho, \theta,x_m) \rangle  \\
      & +\langle \nabla_{x_m} f^i (\rho, \theta, x_m), \tilde{F}_m^i(\rho, \theta,x_m) \rangle\\
      & + \sum_{j \in E} a_{ij}(\rho \theta,x_m) (f^j(\rho, \theta,x_m) - f^i(\rho, \theta,x_m)),
   \end{split}
\end{equation}
where $\tilde{G}^i(\rho, \theta,x_m) = \tilde{F}_n^{i}(\rho,\theta,x_m)-\langle \theta, \tilde{F}_n^{i}(\rho,\theta,x_m) \rangle \theta$.

 The set $\mathcal{X}_0$ is invariant, and we identify it with $S^{n-1} \times \RR^m \times E$. On this set, the process $(\Theta, X^m, I)$ satisfies 

\beq
\label{drho0}
 \left \{
\begin{array}{l}
 \frac{d\Theta_t}{dt} = B^{I_t}(X_t^m) \Theta_t -\langle \Theta_t, B^{I_t}(X_t^m) \Theta_t \rangle \Theta_t \\ \\
  \frac{d X^m_t}{dt} = F_m^{I_t}(0,X^m_t) \\ \\
  \Pr(I_{t+s} = j | {\cal F}_t) = a_{i j}(0, X^m_t) s + o(s) \mbox{ for } i \neq j \mbox{ on } \{I_t = i \}
\end{array}
\right.
\eeq

\blem
\label{lem:extX0}
For all $(\theta,x_m,i) \in \mathcal{X}_0$, one has 
$$\mathbb{P}_{\theta,x_m,i}( \limsup_{t \rar \infty} \frac{1}{t} \log(\|X_t^m\|) \leq \Lambda_D^+) = 1.$$
\elem

\prf On $\mathcal{X}_0$, the process $(X^m,I)$ evolves independently from $\Theta$. It is a PDMP with vector fields $\hat{F}^i : \RR^m \to \RR^m$ and transition rate matrix $(\hat{a}_{ij})$ defined for all $x \in \RR^m$ respectively by $\hat{F}^i(x) = F_m^i(0_n,x)$ and $\hat{a}_{i j}(x) = a_{i j}(0_n, x)$. The origin $0_m$ is a common zero for all the $\hat{F}^i$, and $D \hat{F}^i(0_m) = D^i$. In particular, the maximal Lyapunov exponent for $(X^m,I)$ is $\Lambda_D^+$ and the result follows from  \cite[Theorem 3.1]{BS17} due to the fact that $\Lambda_D^+ < 0$ and $0$ is accessible from $M_0$. \QEDB
\bigskip

Note that on $\{ 0 \} \times S^{n-1} \times \{0_m \} \times E$, $(\Theta,I)$ is equal to the PDMP $(\Theta(B),J)$ defined in section \ref{sec:linear}. Therefore, we have :

\blem
\label{lem:invX0}
Let $\mu$ be an invariant probability of $\tilde{Z}$ on $\mathcal{X}_0$. Then $\mu(\mathrm{d}\theta,\mathrm{d}x,\mathrm{d}i) = \delta_0(\mathrm{d}x)\otimes \hat{\mu}(\mathrm{d}\theta, \mathrm{d}i)$ where $\hat{\mu}$ is an invariant probability of $(\Theta(B),J)$.
\elem 

\prf Let $(Q_t)_{t \geq 0}$ be the semigroup of $(\Theta,X^m,I)$ on $\mathcal{X}_0$. Let $f : \RR^m \to \RR$ be a continuous bounded function and define $\hat{f} : \mathcal{X}_0 \to \RR$ by $\hat{f}(\theta,x,i) = f(x)$. By invariance of $\mu$, $\mu Q_t \hat{f} = \mu \hat{f}$ for all $t \geq 0$. Now, $\mu \hat{f} = \tilde{\mu} f$ where $\tilde{\mu}$ is the marginal of $\mu$ on $\RR^m$ and by Lemma \ref{lem:extX0} and dominated converge, $\mu Q_t \hat{f} \to f(0)$ when $t \to \infty$. Thus $\tilde{\mu}= \delta_0$. Since the marginal law is a Dirac mass, this implies that $\mu$ is a product measure : $\mu = \delta_0 \otimes \hat{\mu}$, where $\hat{\mu}$ is the marginal of $\mu$ on $S^{n-1} \times E$. The result follows from the remark preceding this lemma. \QEDB
\bigskip

Define $H : \mathcal{X} \to \RR$ by $H(\rho,\theta,x_m,i) = - \langle \theta, \tilde{F}_n^i(\rho,\theta,x_m) \rangle$. The following lemma is immediate from Lemma \ref{lem:invX0} and the definition of $H$.

\blem
Let $\mu$ be an invariant probability on $\mathcal{X}_0$. Then with the notation of Lemma \eqref{lem:invX0}, $\mu H = - \Lambda_B( \hat{\mu})$. 
\elem
Now we proceed to the proof of Theorem \ref{thm:pers}. Letting $V : \mathcal{X}_+ \to \RR_+$ be a smooth function coinciding with $- \log(\rho)$ for all $(\rho, \theta, x_m, i) \in \mathcal{X}$ such that $\rho \leq 1$, the end of the proof is verbatim the same as in \cite[Section 5]{BS17} by noting that $\tilde{\mathcal{L}}V = H$ (in a weak sense, see \cite{BS17} or \cite{B18} for details). \QEDB

\subsection{Proof of Proposition \ref{prop:lambdaA}}
The proof is really similar to the one of Theorem \ref{thm:pers}, so we do not give all the details.   Recall from proof of Lemma \ref{lem:triang} that we rewrite \eqref{dThetarho} as

 $$\left \{
\begin{array}{l}
 \frac{d \Theta^n_t}{dt} = B^{J_t} \Theta^n_t -\left( \langle B^{J_t} \Theta^n_t,  \Theta^n_t \rangle +  \langle C^{J_t} \Theta^n_t + D^{J_t} \Theta^m_t, \Theta^m_t\rangle\right) \Theta^n_t \\ \\
 \frac{d \Theta^m_t}{dt} =C^{J_t} \Theta^n_t + D^{J_t} \Theta^m_t -\left( \langle B^{J_t} \Theta^n_t,  \Theta^n_t \rangle +  \langle C^{J_t} \Theta^n_t + D^{J_t} \Theta^m_t , \Theta^m_t \rangle \right) \Theta^m_t
\end{array}
\right.$$
As in the proof of Theorem \ref{thm:pers}, we write $\Theta^n = \rho \hat{\Theta}$ with $\rho = \| \Theta^n \|$ and $\hat{\Theta}=\frac{\Theta^n}{\rho} \in S^{n-1}$. As before, the set $\{ \rho = 0 \}$ is invariant, and one can check that on this state, $\hat{\Theta}$ and $\Theta^m$ evolves independently as : 

\beq
\label{hatteta}
 \left \{
\begin{array}{l}
 \frac{d \hat{\Theta}_t}{dt} = B^{J_t} \hat{\Theta} - \langle B^{J_t} \hat{\Theta},  \hat{\Theta} \rangle  \hat{\Theta}\\ \\
 \frac{d \Theta^m_t}{dt} = D^{J_t} \Theta^m_t - \langle D^{J_t} \Theta^m_t , \Theta^m_t \rangle  \Theta^m_t
\end{array}
\right.
\eeq
That is $(\hat{\Theta},J)=(\Theta(B),J)$ and $(\Theta^m,J)=(\Theta(D),J)$. Furthermore, setting $\hat{V}(\rho,\hat{\theta},\theta^m,i)=-\log(\rho)$, one has 

$$\hat{L}\hat{V}(\rho,\hat{\theta},\theta^m,i)= - \langle B^i \hat{\theta}, \hat{\theta} \rangle + \left(\rho^2 \langle B^i \hat{\theta}, \hat{\theta} \rangle + \rho \langle C^i \hat{\theta}, \theta^m\rangle + \langle D^i \theta^m, \theta^m \rangle \right):= \hat{H}(\rho,\hat{\theta},\theta^m,i).$$  
Here $\hat{L}$ stands for the generator of $\hat{Z}:=(\rho,\hat{\Theta},\Theta^m,J)$. Now if $\mu$ is an invariant probability of $\hat{Z}$ on $\{ \rho = 0 \}$; then there $\hat{\mu} \in \Pcal_{inv}^{(\Theta(B),J)}$ and  $\tilde{\mu} \in \Pcal_{inv}^{(\Theta(D),J)}$ such that $\mu \hat{H} = - \Lambda_B(\hat{\mu}) + \Lambda_D(\tilde{\mu}).$ In particular, if $\Lambda_B^+ < \Lambda_D^-$, then for all $\mu \in \Pcal_{inv}^{\hat{Z}}$ with $\mu(\{ \rho = 0 \}) =1$, one has $\mu \hat{H} > 0$. Moreover, since we assumed that $\{(\theta_n,\theta_m) \in S^{d-1} \: : \: \theta_n=0\}$ is accessible from $S$ for $(\Theta(A),J)$,  the set $\{\rho = 0 \}$ is accessible for $\hat{Z}$.  This concludes the proof by the same arguments as in \cite{BS17}. 

\brem
{\rm The same proof shows that if $\Lambda_B^- > \Lambda_D^+$, $(\Theta(A),J)$ admits at least one invariant probability measure giving no mass to $\{(\theta_n,\theta_m) \in S^{d-1} \: : \: \theta_n=0\}$. An interesting question would be to know if it is possible to recover the Lyapunov exponents associated to $B$ with this invariant probability measure, like in dimension $2$ (see Example \ref{ex:dim2}.) }
\erem

\subsection{Proof of Proposition \ref{prop:triang}}
For this proof, we use a result from Hennion \cite{H84} which is given for Lyapunov exponents of random dynamical systems. Recall from \cite{BS17} that if $\Omega = \mathbb{D}(\RR_+,E)$, the set of cadlag functions with values in $E$, then $(\Omega, \mathcal{F}, (\mathbf{\Theta})_{t \geq 0}, \PP_p^J)$ is an ergodic dynamical system. Here, $\mathcal{F}$ is the Borel sigma field of $\Omega$ and $\mathbf{\Theta}_t : \Omega \to \Omega$ is the $t$-time shift defined for all $\omega \in \Omega$ by $\mathbf{\Theta}_t(\omega_s)= \omega_{t+s}$.For a set of matrices $(\hat{A}^i)_{i \in E}$ and for $\omega \in \Omega$ and $y \in \RR^d,$ let $$t \mapsto \varphi(t,\omega) y$$ denote  the solution to the linear  differential equation
$$\dot{y} = \hat{A}^{\omega_t} y$$ with initial condition $\varphi(0,\omega) y = y.$

Then,
$(\varphi,(\mathbf{\Theta})_{t \geq 0})$ is a {\em linear random dynamical system} over the ergodic dynamical system
$(\Omega, \mathcal{F}, (\mathbf{\Theta})_{t \geq 0}, \PP_p^J)$. Furthermore, the conditions of the Multiplicative Ergodic Theorem are satisfied (see \cite{ColMa15} for details) :

\bthm[Multiplicative Ergodic Theorem]\

There exist $1 \leq \tilde{d} \leq d,$ numbers
$$\lambda_{\tilde{d}} < \ldots < \lambda_1,$$ called {\em the Lyapunov exponents} of $(\varphi, \mathbf{\Theta})$ , a Borel set $\tilde{\Omega}  \subset \Omega$ with $\mathbb{P}^J_{p}(\tilde{\Omega}) = 1,$ and for each $\omega \in \tilde{\Omega}$ distinct vector spaces
$$\{0\} = V_{\tilde{d}+1}(\omega) \subset V_{\tilde{d}}(\omega) \subset \ldots \subset V_{i}(\omega) \ldots \subset V_1(\omega) = \RR^d$$ (measurable in $\omega$) such that
\beq
\label{eq:deflambdai}
\lim_{t \rar \infty} \frac{1}{t} \log \|\varphi(t,\omega) y\| = \lambda_i
\eeq for all $y \in V_i(\omega) \setminus V_{i+1}(\omega)$.
\ethm
We let $S(\hat{A})$ denoted the set of Lyapunov exponents given by the above theorem and $\lambda_1(\hat{A}) = \max S(\hat{A})$. By \cite[Proposition 2.5]{BS17},  $\Lambda(\hat{A}) \subset \mathcal{S}(\hat{A})$ and $\Lambda_{\hat{A}}^+=\lambda_1(\hat{A})$. Now  \cite[Proposition 1]{H84} implies that $\mathcal{S}(A)=\mathcal{S}(B)\cup \mathcal{S}(D)$ and thus $\lambda_1(A)=\max( \lambda_1(B), \lambda_1(D) )$. Hence the result. \QEDB

\brem
{\rm Proposition 1 in \cite{H84} is given for a discrete-time random dynamical system. However, its proof adapts verbatim to the continuous-time case by using the continuous-time version of the Multiplicative Ergodic Theorem \cite[Theorem 3.4.1 (C)]{Arn98}. One could also have used \cite[Theorem 1.1]{GMO08}.}
\erem

\section*{Acknowledgments}
This work was supported by the SNF grant $2000020 - 149871/1$. I warmly thank Michel Bena{\"i}m and Tobias Hurth for helpful discussions and review of the paper.

\bibliographystyle{amsalpha}
\bibliography{RandomSwitch}
\end{document}